\documentclass[draft]{article}
\usepackage{amssymb}
\usepackage{breakcites}
\usepackage{amsmath}
\usepackage{amsthm}
\usepackage{amsfonts}
\usepackage{mathrsfs}
\usepackage{amssymb,color,dsfont}
\allowdisplaybreaks
\usepackage{bm}

\usepackage{amsfonts, txfonts, graphicx, nicefrac}

    \def\fz{\varphi}
  \def\kz{\kappa}
\def\lz{\lambda} \def\mz{\mu}
\def\nz{\nu}     
\def\pz{\pi}

\def\rz{\rho}        
        
\def\vz{\varepsilon}

\def\llz{\Lambda} \def\ddz{\Delta}

\def\aq{{\mathscr{A}}}

   \def\pq{{\mathscr{P}}}

\def\qd{\quad}
\def\qqd{\qquad}

\def\lt{\left}
\def\rt{\right}

\def\PP{\mathbb{P}}
\def\EE{\mathbb{E}}

\newcommand{\mathsym}[1]{{}}

\def\le{\leqslant}
\def\ge{\geqslant}
\def\leq{\leqslant}
\def\geq{\geqslant}

\newtheorem{thm}{Theorem}[section]
\newtheorem{prop}[thm]{Proposition}

\newtheorem{rem}[thm]{Remark}

\newtheorem{conj}[thm]{Conjecture}

\numberwithin{equation}{section} \allowdisplaybreaks[4]

\def\deprf{\quad $\square$ \medskip}
\def\bg{\begin}
\def\be{\begin{equation}}
\def\de{\end{equation}}
\def\ben{\begin{equation*}}
\def\den{\end{equation*}}

\def\dear{\end{eqnarray*}}
\def\lb{\label}

\def\dps{\displaystyle}

\def\ben{\bg{enumerate}}
\def\den{\end{enumerate}}

\def\d{\text{\rm d}}

\def\law{\mathcal{L}}

\def\d{\text{\rm d}}

\def\PP{\mathbb{P}}
\def\EE{\mathbb{E}}
\def\RR{\mathbb{R}}
\def\ZZ{\mathbb{Z}_+}

\def\LF{\overleftarrow{F}}
\def\RF{\overrightarrow{F}}

\def\Lip{\text{\rm Lip} (\rz)}

\def\Ref#1{(\ref{#1})}

\def\sW{\mathscr{W}}
\def\mW{\mathbb{W}}
\def\sLo{\mathscr{L}^0}
\def\bone{{\bf 1}}

\begin{document}
\date{}
\pagestyle{plain}
\title{On Stein's factors for Poisson approximation in Wasserstein distance with non-linear transportation costs}
\author{Zhong-Wei Liao \footnote{Postal address: South China Research Center for Applied Mathematics and Interdisciplinary Studies, South China Normal University, Guangzhou 510631, China. (zhwliao@m.scnu.edu.cn)}, Yutao Ma\footnote{Postal address: School of Mathematical Sciences, Beijing Normal University, Beijing 100875, China. (mayt@bnu.edu.cn)}, Aihua Xia \footnote{Postal address: School of Mathematics and Statistics, The University of Melbourne, VIC 3010, Australia. (aihuaxia@unimelb.edu.au)}}
\date{}
\maketitle \underline{}

{\bf Abstract:} We establish various bounds on the solutions to a Stein equation for Poisson approximation in Wasserstein distance with non-linear transportation costs.
The proofs are a refinement of those in \cite{BX06} using the results in \cite{LM09}. As a corollary, we obtain an estimate of Poisson approximation error measured in $L^2$-Wasserstein distance.



\setlength{\baselineskip}{0.25in}

\section{Framework and introduction} \label{sect1}

As the cornerstone of the law of small numbers, Poisson distribution provides good approximation to the distribution of the counts of rare events and the quality of Poisson approximation has been studied extensively in the literature~\cite{BHJ92}. In particular, the pioneering works of \cite{C75,B88} enable us to assess the accuracy of Poisson approximation to the distribution of the sum of integer valued random variables under a variety of dependent structures in terms of various metrics. The key to the success is the so called Stein's factors. When the approximation errors are measured in the total variation distance, \cite{BH84} conclude that sharp bounds of Stein's factors often yield remarkably sharp estimates of the approximation errors. However, sharp estimates of Stein's factors for Poisson approximation are generally hard to extract and, in addition to the total variation distance and the Kolmogorov distance, the only conclusive case is in terms of the Wasserstein distance with linear transportation costs~\cite{BX06}. In the field of mass transportation problems, the Wasserstein distance plays a pivotal role but the transportation costs are often non-linear~\cite{V03}. For example, what is the {$L^2$}-Wasserstein distance between a Poisson binomial distribution and a Poisson distribution? In this paper, we aim to tackle the  problem and establish various bounds on the solutions to a Stein equation for Poisson approximation in terms of the Wasserstein distance with non-linear transportation costs. The bounds are used to quantify the accuracy of Poisson approximation to the Poisson binomial distribution in $L^2$-Wasserstein distance.

Given any $\lz > 0$, denote by $\pz_i = e^{-\lz} \lz^i / i!$, $i\in\ZZ:=\{0,1,2,\dots\}$, the Poisson distribution with mean $\lz$. Denote by $\pq (\ZZ)$ the set of all probability measures on $\ZZ$ and $\aq$ the set of all strictly increasing functions $\rz$ on $\ZZ$ such that $\sum_{i=0}^\infty|\rz(i)|\pz_i<\infty$. Each $\rz \in \aq$ induces a metric on $\ZZ$ through
$$
d_{\rz} (i, j) = |\rz(i) - \rz(j)|, \qqd \forall i, j \in \ZZ.
$$
The Wasserstein distance between $\nz_1,\nz_2 \in \pq (\ZZ)$ with non-linear transportation costs considered in the paper is defined by
$$
\sW_{d_{\rz}} (\nz_1, \nz_2) = \inf  \sum_{i, j \in \ZZ} d_{\rz} (i, j) \mz (i, j),
$$
where the infimum is taken over all couplings $\mz$ of $\nz_1$ and $\nz_2$ such that  $\nz_1(\cdot )= \mz(\cdot, \ZZ)$ and $\nz_2(\cdot) = \mz (\ZZ, \cdot)$. Obviously, when $\rz (i) = i$, the distance $\sW_{d_{\rz}}$ degenerates to $L^1$-Wasserstein distance, i.e., with linear transportation costs. The Kantorovich-Rubinstein duality theorem \cite{KR58,E11} says that
\be\lb{dual}
\sW_{d_{\rz}} (\nz_1, \nz_2) = \sup \bigg\{ \nz_1(f) - \nz_2 (f) : \| f \|_{\Lip} = 1 \bigg\},
\de
where $\nu_j (f) := \sum_{i \in \ZZ} f(i) \nu_j(\{i\})$ for $j=1,2$ and
$$
\| f \|_{\Lip} := \sup_{i \neq j} \frac{|f (j) - f (i)|}{|\rz (j) - \rz(i)|} = \sup_{i \geq 0} \frac{|f (i+1) - f (i)|}{\rz (i+1) - \rz(i)}.
$$
A function $f$ on $\ZZ$ is called {\it $\rz$-Lipschitzian} if $\| f \|_{\Lip} < \infty$ and one can easily verify that $\| f \|_{\Lip} = 1$ in \Ref{dual} can be replaced with $|f(i) - f(j)| \leq |\rz(i) - \rz(j)|$, $\forall i, j \in \ZZ$. The duality form (\ref{dual}) has a long history, dating back to \cite{KR58} on the mass transport problems, see \cite[Chapter 5]{RKSF13} for more details. The metric $\sW_{d_{\rz}}$ belongs to the family of the $L^1$-Wasserstein distance and it remains an open problem to use Stein's method for estimating approximation errors in terms of other $L^p$-Wasserstein distances ($1< p < \infty$) for probability measures $\nu_1$ and $\nu_2$ on $\RR$ defined by
$$
\mW_p(\nu_1,\nu_2)=\left(\inf \int_{{\RR \times \RR}} |x-y|^p\mz(dx,dy)\right)^{1/p},
$$
where, as before, the infimum is taken over all couplings $\mz$ of $\nz_1$ and $\nz_2$ with  $\nz_1(\cdot )= \mz(\cdot, \RR)$ and $\nz_2(\cdot) = \mz (\RR, \cdot)$. This is because the Kantorovich-Rubinstein duality theorem for $\mW_p$ with $p\ne 1$ does not possess the form \Ref{dual} which is the key to the Stein equation~\Ref{stein}. Nevertheless, since $(i-j)^2\le \left|i^2-j^2\right|$ for all $i,j\in\ZZ$, we have the following crude estimate for $\mW_2$.

\begin{prop} \label{W2} For any two probability measures $\nu_1,\nu_2$ on $\ZZ$, with $\rz_2(\cdot)=\cdot^2$, we have
$$\mW_2(\nu_1,\nu_2)\le \left(\sW_{d_{\rz_2}} (\nz_1, \nz_2) \right)^{1/2}.$$
\end{prop}

For any random variable $W$ on $\ZZ$, the Stein-Chen method for estimating the distance between the distribution $\law(W)$ of $W$ and $\pi$ is based on the following observation~\cite{C75}: $W$ follows the distribution $\pz$ if and only if
\begin{equation}
\EE \lt[ \lz g (W+1) - W g(W) \rt] =0,\label{Steinop}
\end{equation}
for all functions $g: \ZZ \to \RR$ satisfying $\EE \lt[ W |g| (W)\rt] < \infty$. This leads to the well-known Stein equation for Poisson approximation: for each $f$ on $\ZZ$,
\be\lb{stein}
\lz g_f (i+1) - i g_f (i) = f (i) - \pz(f), \qqd i \in \ZZ,
\de
and one can recursively solve for the function $g_f$. As the value of $g_f(0)$ does not affect the equation, we set $g_f (0) := g_f(1)$ for convenience. Using (\ref{dual}) and (\ref{stein}), the $\sW_{d_{\rz}}$ distance between $\law (W)$ and $\pz$ can be reformulated as
\be\lb{stein1}
\sW_{d_{\rz}} (\law (W), \pz)  = \sup_{\| f \|_{\Lip} = 1} | \EE [ f (W)] - \pz (f) | = \sup_{\| f \|_{\Lip} = 1} \lt| \EE \lt[ \lz g_f (W+1) - W g_f (W) \rt] \rt|.
\de
On the other hand, one can often use the dependence structure of $W$ to expand the right-hand side of (\ref{stein1}) into
$$
\lt| \EE \lt[ \lz g_f (W+1) - W g_f (W) \rt] \rt| \leq \EE\vz_0 M_0 (g_f) + \EE\vz_1 M_1 (g_f) + \EE\vz_2 M_2 (g_f),
$$
where $\EE\vz_k\ge 0$,
 \begin{equation}
{M_k (g_f) = {\sup_{i \geq 1}} \frac{|\ddz^k g_f (i)|}{\ddz\rz(i)}}, \qqd k=0,1, 2, \label{def:Mk}
\end{equation}
and $\ddz$ is the difference operator defined as $\ddz g (i) = g(i+1) - g(i)$ and $\ddz^k g(i) = \ddz^{k-1} g (i+1) - \ddz^{k-1} g (i)$, $k \geq 2$.
This, together with \Ref{stein1}, ensures
$$
\sW_{d_{\rz}} (\law (W), \pz) \leq \EE\vz_0 \sup_{\| f \|_{\Lip} = 1} M_0 (g_f) +\EE \vz_1 \sup_{\| f \|_{\Lip} = 1} M_1 (g_f) + \EE\vz_2 \sup_{\| f \|_{\Lip} = 1} M_2 (g_f).
$$

The birth-death process interpretation of $g_f$ in \cite{B88} says if we write $g_f (i) = h_f (i) - h_f (i-1)$, then Stein's equation (\ref{stein}) becomes
\be\lb{stein2}
\lz (h_f (i+1) - h_f (i)) - i (h_f (i) - h_f (i-1)) = f(i) - \pz(f), \qqd \forall i \geq 1.
\de
This ensures that $h_f$ is the solution to the Stein equation (which is also known as Poisson equation)
\be\lb{Poisson}
Q h_f = f - \pz(f),
\de
where $Q$ is a transition matrix defined as
\begin{align*}
q_{i, i+1} & = \lz > 0, \qd q_{i,i}= - (\lz +i), \qd \forall i \geq 0; \qd q_{i, i-1} = i, \qd \forall i \geq 1, \\
q_{i, j} &= 0, \qd \text{if $|i-j|>1$, for $i, j \in \ZZ$.}
\end{align*}
Denote by $\sLo (\rz)$ the space of $\rz$-Lipschitzian functions $f$ satisfying $\pz (f) = 0$. The definition of $Q$ ensures that the unique solution to the equation $Q h = 0$ with $\pz (h) =0$ is $h\equiv 0$. Hence, for each $f \in \sLo (\rz)$, there exists a unique solution $h_f$ with $\pz (h_f) =0$ to the equation $Q h_f = f$, which means that $Q^{-1}$ is well defined on $\sLo (\rz)$. Moreover, the operator norm of $(-Q)^{-1}$ is defined as
$$
\| (-Q)^{-1} \|_{\Lip} := \sup \lt\{ \| (-Q)^{-1} (f - \pz(f)) \|_{\Lip}: \  \| f \|_{\Lip} =1 \rt\}.
$$
See \cite{C10} and \cite{LM09} for more information of the Poisson equation and the spectral gap of birth-death processes.

The upper bounds of Stein's factors $\sup_{\| f \|_{\Lip} = 1} M_k (g_f)$ for $\sW_{d_{\rz}}$ distance are summarized in the following theorem.

\begin{thm}\lb{thmPo}
{
{Let $\rz \in \aq$, $h_{\rz}$ be the solution to equation $Q h_{\rz} = \rz - \pz(\rz)$ and $\lfloor \lz \rfloor$ be the largest integer less than or equal to $\lz$. Define $m_\rz=\sup_{i\ge 0}\frac{\ddz \rz(i)}{\ddz\rz(i+1)}$}.  Then we have
\begin{align}
\sup_{\| f \|_{\Lip} = 1} M_0 (g_f) & \leq m_\rz\lt\| (-Q)^{-1} \rt\|_{\Lip}, \lb{Po0}\\
\sup_{\| f \|_{\Lip} = 1} M_2 (g_f) & \leq m_\rz\lt\| \ddz^2 h_{\rz} \rt\|_{\Lip} + 2 \lt( (2 \Xi_2 (\lz)) \wedge \lz^{-1} \rt), \lb{Po2}
\end{align}
where
\begin{align}
\Xi_2 (\lz) & :=
\left\{
  \begin{array}{ll}
    \dps\frac{(\lz -1)^2 -2 e^{-\lz} +1}{\lz^3}, & 0 < \lz \leq 1, \\[+5pt]
    \dps\frac{(e - 1) (\lz-1)^2 + 2 \lz + e -4}{\lz^3 e} + \dps\sum_{n=1}^{\lfloor \lz \rfloor -1} \frac{4 \sqrt{n} (3 (\lz - n)^2 - 3(\lz - n) + 1)}{ \sqrt{2 \pz} \lz^3 (12n + 1)} & \\
  \qqd + \dps\frac{4 \sqrt{\lfloor \lz \rfloor} ( \lz - \lfloor \lz \rfloor)^3}{\sqrt{2 \pz} \lz^3 (12 \lfloor \lz \rfloor + 1)}, & 1< \lz < \infty,
  \end{array}
\right. \lb{lz2}\\
& \le
\left\{
  \begin{array}{ll}
    \dps\frac{1}{3}, & 0 < \lz \leq 1, \\[+5pt]
    \dps\frac{0.426}{\sqrt{\lz}}, & 1< \lz < \infty.
  \end{array}
\right. \lb{xialz2}
\end{align}
If $\ddz^2 \rz (i) \geq 0$, $\forall i \in \ZZ$, then
\begin{align}
\sup_{\| f \|_{\Lip} = 1} M_1 (g_f)  & \leq m_\rz\lt\| \ddz h_{\rz} \rt\|_{\Lip} + 2m_\rz \Xi_1 (\lz), \lb{Po1}
\end{align}
and if $\ddz^2 \rz (i) \leq 0$, $\forall i \in \ZZ$, then
\begin{align}
\sup_{\| f \|_{\Lip} = 1} M_1 (g_f)  & \leq m_\rz\lt\| \ddz h_{\rz} \rt\|_{\Lip} + 2 \Xi_1 (\lz), \lb{XiaPo1}
\end{align}
where
\begin{align}
\Xi_1 (\lz) & :=
\left\{
  \begin{array}{ll}
    \dps\frac{e^{- \lz} + \lz - 1}{\lz^2}, & 0 < \lz \leq 1, \\[+5pt]
    \dps\frac{(e - 1) (\lz -1) +1 }{\lz^2 e} + \dps\frac{6 \sqrt{ \lfloor \lz \rfloor} (\lz - \lfloor \lz \rfloor)^2 }{\sqrt{2 \pz} \lz^2 (12 \lfloor \lz \rfloor + 1)} + \dps\sum_{n=1}^{\lfloor \lz \rfloor -1} \frac{ 12 \sqrt{n} (\lz -n) - 6 \sqrt{n} }{ \sqrt{2 \pz} \lz^2 (12 n + 1)}, & 1< \lz < \infty,
  \end{array}
\right. \lb{lz1} \\
\notag \\
& \le
\left\{
  \begin{array}{ll}
    \dps\frac{1}{2}, & 0 < \lz \leq 1, \\[+5pt]
    \dps\frac{0.532}{\sqrt{\lz}}, & 1< \lz < \infty.
  \end{array}
\right. \lb{xialz1}
\end{align}
}
\end{thm}

\begin{rem}\lb{rem0}
According to \cite[Lemma 2.3]{LM09}, $\ddz h_{\rz}$ mentioned above is explicit and computable,
\be\lb{hrz}
\ddz h_{\rz} (i) = h_{\rz} (i+1) - h_{\rz} (i) = \frac{1}{(i+1) \pz_{i+1}} \sum_{j=0}^{i} \pz_j (\rz(j) - \pz (\rz)), \qqd i \geq 0.
\de
Moreover, $h_{\rz}$ {has a simple and straightforward expression for many cases}, see Proposition \ref{prop} below.
\end{rem}

Recalling the definition of $M_k (g_f)$ in \Ref{def:Mk}, we can see that $\ddz^kg_f(0)$ is excluded in the definition. This is because the value of $g_f (0)$ has no effect on the Stein equation \Ref{stein} and we can set it to any value. However, whatever value we set for $g_f(0)$, there is a direct consequence on $\ddz^kg_f(0)$ for $k\ge 0$ and there seems to be no optimal values such that we can incorporate them into the bounds in Theorem~\ref{thmPo}. Here we consider the approach in \cite{BX06} with the following bounds.

\begin{prop}\lb{remadd} With $g_f (0) = g_f (1)$, we have
\begin{align}
\sup_{\| f\|_{\Lip} = 1} \frac{|g_f (0)|}{\ddz \rz (0)} & = \frac{\pz (\rz) - \rz (0)}{\lz \ddz \rz(0)}, \lb{Mg0} \\
\sup_{\| f\|_{\Lip} = 1} \frac{|\ddz g_f (0)|}{\ddz \rz (0)} & = 0, \lb{Mg1} \\
\sup_{\| f\|_{\Lip} = 1} \frac{|\ddz^2 g_f (0)|}{\ddz \rz (0)} & = \lt| \frac{1}{\lz} + \frac{\rz(0) - \pz (\rz)}{\lz^2 \ddz \rz (0)} \rt| +
\left\{
  \begin{array}{ll}
    \dps \frac{2(e^{-\lz} + \lz -1)}{\lz^2}, & \hbox{when $\ddz^2 \rz (\cdot) \geq 0$}; \\
    \dps \frac{2 \ddz \rz (1) (e^{-\lz} + \lz -1)}{\ddz \rz (0) \lz^2}, & \hbox{when $\ddz^2 \rz (\cdot) \leq 0$.}
  \end{array}
\right. \lb{Mg2}
\end{align}
\end{prop}

\begin{rem}\lb{rem1}
We can directly verify that $\rz_1 (i) = i $ satisfies
\be\lb{Poisson1}
-Q \rz_1 = \rz_1 - \pz(\rz_1), \qd \text{\rm and} \qd \rz_1 \in \aq,
\de
which implies that $\rz_1 - \pz(\rz_1)$ is the eigenfunction of $-Q$ corresponding to the eigenvalue $\kz = 1$. By \cite[Theorem 3.1]{LM09}, $\lt\| (-Q)^{-1} \rt\|_{\Lip}$ attains the supremum at the eigenfunction of $-Q$ and equals to the reciprocal of eigenvalue $\kz^{-1} = 1$. In this cae, $m_{\rz_1}=1$, the distance $\sW_{d_{\rz}}$ is consistent with the $L^1$-Wasserstein distance studied in \cite{BX06}, and the bounds (\ref{Po0}) is the same as the result given in \cite[Theorem 1.1]{BX06}.
\end{rem}

\begin{rem}\lb{rem2}
When $\rz = \rz_1$, by (\ref{Poisson1}), we have $h_{\rz_1} = - i$ and then $\lt\| \ddz h_{\rz_1} \rt\|_{\text{\rm Lip}(\rz_1)} = 0$. Hence,
\be\lb{M1}
\sup_{\| f \|_{{\rm Lip} (\rz_1)} = 1} M_1 (g_f) \le
\left\{
  \begin{array}{ll}
    (e^{-\lz} + \lz -1) / \lz^2, & \hbox{for $0<\lz \le 1$,} \\[+5pt]
    \dps\frac{1.064}{\sqrt{\lz}}, & \hbox{for $1<\lz < \infty$.}
  \end{array}
\right.
\de
It should be pointed out that when $0<\lz \leq 1$ the estimate of $\sup_{\| f \|_{{\rm Lip} (\rz_1)} = 1} M_1 (g_f)$ is sharp (see (\ref{Po1add2}) below), and when $1<\lz < \infty$ the constant of the estimate slightly improves \cite[Theorem 1.1]{BX06}. The function $\Xi_2 (\lz)$ has the same order as that of $\Xi_1 (\lz)$ for $\lz\to\infty$. When $\rz =\rz_1$, we have $\lt\| \ddz^2 h_{\rz_1} \rt\|_{\text{\rm Lip}(\rz_1)} = 0$ and
\be\lb{M2}
\sup_{\| f \|_{{\rm Lip} (\rz_1)} = 1} M_2 (g_f) \le
\left\{
  \begin{array}{ll}
    4 \lt( (\lz-1)^2 - 2 e^{-\lz} + 1 \rt) / \lz^3, & \hbox{for $0<\lz \le1$,} \\[+5pt]
    \dps\frac{1.704}{\sqrt{\lz}}\wedge \frac2 \lz, & \hbox{for $\lz >1$,}
  \end{array}
\right.
\de
hence (\ref{M2}) is slightly better than \cite[Theorem~1.1]{BX06} but with the same asymptotic behaviour when $\lz$ is close to $0$ or is large.
\end{rem}

The Wasserstein distance in Theorem~\ref{thmPo} covers a range of cost functions and one can choose different $\rz$ depending on the problem of interest. We demonstrate how to solve \Ref{hrz} in the following proposition.

\begin{prop}\lb{prop}
{
(1) Consider the convex case $\rz_p (i) := i^p$, where $p \geq 1$. Denote by $h_p$ the solution to the Stein equation $Q h_p =\rz_p - \pz(\rz_p)$. Then for each $i \geq 1$, $h_p (i)$ satisfies the recursive formula
\be\lb{ffz}
h_p (i) =
\left\{
  \begin{array}{ll}
    -i, & \hbox{$p=1$;} \\[+5pt]
    - \dps\frac{i^p}{p} + \dps\frac{1}{p} \dps\sum_{k=1}^{p-1} \binom{p}{k} h_k (i) \lt[ \lz + \dps\frac{k (-1)^{p-k+1}}{p-k+1} \rt], & \hbox{$p \geq 2$,}
  \end{array}
\right.
\de
and $h_p (0) = h_p (1) + \lz^{-1} \pz (\rz_p)$. In particular, when $p=2$, it implies that for each $i \in \ZZ$, $\ddz \rz_{2} (i) \geq 0$, $\ddz^2 \rz_{2} (i) \geq 0$ and $m_{\rz_2}=1$, giving
\be\lb{rz2}
\sup_{ \| f \|_{{\rm Lip} (\rz_2)} = 1} M_0 (g_f) \leq \lz +1, \qd \sup_{\| f \|_{{\rm Lip} (\rz_2)} = 1} M_1 (g_f) \leq 1 + 2 \Xi_1 (\lz), \qd \sup_{\| f \|_{{\rm Lip} (\rz_2)} = 1} M_2 (g_f) \leq {2 \lt( (2\Xi_2 (\lz)) \wedge \lz^{-1} \rt)}.
\de

(2) Consider the concave case $\rz_{1/2} (i) := {\lz +} \sqrt{i} - \lz / \sqrt{i+1}$, it implies that for each $i \geq 0$, $\ddz \rz_{1/2} (i) \geq 0$, $\ddz^2 \rz_{1/2} (i) \leq 0$ and
{$$
m_{\rz} = \frac{\sqrt{3} \lt( \sqrt{2} + \sqrt{2}\lz - \lz \rt)}{\sqrt{3} \lt(2-\sqrt{2} \rt) + \lz \lt( \sqrt{3}-\sqrt{2} \rt)}.
$$}
Then,
\begin{align*}
& \sup_{\| f \|_{{\rm Lip} (\rz_{1/2})} = 1} M_0 (g_f) \leq 2 {m_{\rz}}, \\
& \sup_{\| f \|_{{\rm Lip} (\rz_{1/2})} = 1} M_1 (g_f) \leq  \frac{{m_{\rz}}}{\lz + \lt(\sqrt{2} + \sqrt{3} \rt) \lt(2 \sqrt{3} - \sqrt{6} \rt)} + 2 \Xi_1 (\lz), \\
& \sup_{\| f \|_{{\rm Lip} (\rz_{1/2})} = 1} M_2 (g_f) \leq \frac{\lt(\sqrt{2} + 1 \rt) \lt(2+ \sqrt{2} - \sqrt{6}/3 \rt)}{\lz +2 + \sqrt{2}} {m_{\rz}} + 2 \lt( (2\Xi_2 (\lz)) \wedge \lz^{-1} \rt).
\end{align*}}
\end{prop}

As in \cite{BX06}, we use Poisson approximation to the Poisson binomial distribution to show the accuracy of the bounds for $\rz_2(i)=i^2$.

\begin{prop} \label{example} Let $X_i,\ 1\le i\le n, $ be independent Bernoulli random variables with $\EE X_i=p_i$ and define $W = \sum_{i=1}^n X_i$,
$\mu= \sum_{i=1}^n p_i$, $\mu_l := \sum_{i=1}^n p_i^l$, $\lambda= \mu- \mu_2$.  If $\mu_2$ is
an integer, then we have 
\begin{equation}
\sW_{d_{\rz_2}} (\law ((W-\mu_2)\bone_{W\ge \mu_2}), \pz) \le {6}(\mu_2-\mu_3)+\mu_2({7}+\lambda)e^{-\lambda^2/(2\mu)} \label{example1}
\end{equation}
and
\begin{equation}
\mW_2(\law (W), \pz*\delta_{\mu_2})=\mW_2(\law (W-\mu_2), \pz)\le \mu_2e^{-\lambda^2/(4\mu)}+ \left\{{6}(\mu_2-\mu_3)+\mu_2({7}+\lambda)e^{-\lambda^2/(2\mu)}\right\}^{1/2},  \label{example2}
\end{equation}
where $\delta_{\mu_2}$ is the Dirac measure at $\mu_2$ and $*$ denotes convolution.
\end{prop}

\begin{conj} We conjecture that the order of the upper bound in \Ref{example2} can be significantly improved.
\end{conj}

\section{The proofs}

We first note that \Ref{xialz2} and \Ref{xialz1} are obtained from a numerical computation. For the remaining claims, we need the following notations and preliminaries. Denote by $(X_t^i)_{t \geq 0}$ the birth-death process corresponding to $Q$ with the initial value $X_0^i = i$. Let $P_t$ be the semigroup of $X_t^i$. By \cite{BX06} or \cite{BX01}, we can couple $X_t^{i}$ and $X_t^{i-1}$ by setting
\be\lb{couple_Po}
X_t^{i} = X_t^{i-1} + \mathbf{1}_{ \{\llz > t \} }, \qqd t \geq 0, \qd i \geq 1,
\de
where $\llz$ is a negative exponential random variable with mean $\EE [\llz] = 1$ and independent of $X_t^{i-1}$. According to \cite[Chapter 3.2]{A91}, for any $i \in \ZZ$, we have the expression of the semigroup of $X_t^i$
\be\lb{PoPt}
P_t (i, j) = e^{-\lz (1- e^{-t})} \sum_{k=0}^{i \wedge j} \frac{i !}{k! (i-k)! (j-k)!} e^{-k t} (1- e^{-t})^{i-k} (\lz - \lz e^{-t})^{j-k}, \qqd t \geq 0, \ i, j \in \ZZ.
\de
By integration by parts, it is easy to verify that
\be\lb{I-Q}
\int_0^\infty e^{-t} P_t f (i) \d t = (I - Q)^{-1} f (i), \qd \forall i \in \ZZ,
\de
whenever the integral is well-defined. Moreover, using \Ref{stein2}, we have
\begin{align*}
f (i+1) - f (i)
&= - \ddz h_{f} (i) + \lz (\ddz h_{f} (i+1) - \ddz h_{f} (i)) + i (\ddz h_{f} (i-1) - \ddz h_{f} (i)) \\
&= - \ddz h_{f} (i) + Q (\ddz h_{f} ) (i) = - (I-Q) (\ddz h_{f} ) (i),
\end{align*}
giving
\be\lb{I-Q2}
\ddz h_{f}(i)=- (I-Q)^{-1}(\ddz f)(i), \qqd i \in \ZZ.
\de

Denote
\be\lb{e0}
e_i^+ = (\lz \pz_i)^{-1} \LF (i), \qqd e_i^- = (i \pz_i)^{-1} \RF (i),
\de
where $\LF (i) = \dps\sum_{k=0}^i \pz_k$ and $\RF (i) = \dps\sum_{k=i}^\infty \pz_k$. According to \cite[Lemma 2.4]{BX01} and \cite[Lemma 2.1 and p.~950]{BX06}, for each $i \geq 1$, we have
\begin{align}
\ddz e_i^+ & := e_{i+1}^+ - e_i^+ \geq 0, \qd \ddz e_i^- := e_{i+1}^- - e_{i}^- \leq 0, \lb{e1} \\
\ddz^2 e_{i-1}^+ & := e_{i+1}^+ - 2 e_i^+ + e_{i-1}^+ \geq 0, \qd \ddz^2 e_i^- := e_{i+2}^- - 2e_{i+1}^- + e_i^- \geq 0, \lb{e2}\\
r_i&:=\pz_{i+1}(2e_i^+-e_{i-1}^++e_{i+2}^-)-(e_{i+1}^+-2e_i^++e_{i-1}^+)\RF (i+2) \notag \\
&\ = (- \ddz^2 e_{i-1}^+) \RF(i+1) + \lz^{-1} \geq0.\lb{xiari}
\end{align}
Having these in mind, we are ready to prove the main theorem.

\vskip 1 cm

\noindent {\bf Proof of (\ref{Po0}).} Since $h_f$ is the solution to the Stein equation (\ref{Poisson}), using \cite[Lemma 2.3]{LM09}, we have
\be\lb{Po0-0}
g_f (i) = h_f (i) - h_f (i-1) = \frac{1}{i \pz_i} \sum_{j=0}^{i-1} \pz_j (f(j) - \pz (f)), \qd {i \geq 1}.
\de
On the other hand, according to \cite{B88} or \cite{BX01}, $h_f$ can be expressed as
\be\lb{Po0-1}
h_f (i) = - \int_0^\infty \lt[ \EE [f (X_t^i)] - \pz (f) \rt] \d t.
\de
By the coupling in (\ref{couple_Po}), we have from (\ref{Po0-1}) that
$$
g_f (i) =- \int_0^\infty \left\{\EE [f (X_t^i)] - \EE [f (X_t^{i-1})] \right\}\d t = - \int_0^\infty e^{-t} \EE \lt[ f(X_t^{i-1} + 1) - f(X_t^{i-1}) \rt] \d t,
$$
which implies that
$$
\sup_{\| f \|_{{\rm Lip} (\rz)} = 1} |g_f (i)| = \int_0^\infty e^{-t} \EE \lt[ \ddz \rz (X_t^{i-1}) \rt]dt,
$$
where the supremum is attained by $f = - \rz$. Hence, by (\ref{Po0-0}) we have
$$
\sup_{\| f \|_{{\rm Lip} (\rz)} = 1} |g_f (i)| = g_{-\rz} (i) = \ddz h_{-\rz} (i-1) = \frac{1}{i \pz_i} \sum_{j=0}^{i-1} \pz_j (- \rz(j) + \pz (\rz)) = \frac{1}{i \pz_i} \sum_{j=i}^{\infty} \pz_j (\rz(j) - \pz (\rz)).
$$
Using the representation of $\| (-Q)^{-1} \|_{\Lip}$ given in \cite[Theorem 2.1]{LM09}, it holds that
{
\begin{align}\lb{Po0-2}
\| (-Q)^{-1} \|_{\Lip} &= \sup_{i \geq 1} \frac{\sum_{j=i}^{\infty} \pz_j (\rz (j) - \pz (\rz))}{i \pz_i (\rz(i)- \rz (i-1))} = \sup_{i \geq 1} \frac{\ddz h_{-\rz} (i-1)}{\ddz \rz(i)} \cdot \frac{\ddz \rz(i)}{\ddz \rz (i-1)} \notag \\
&\geq \lt( \sup_{i \geq 1} \sup_{\| f \|_{{\rm Lip} (\rz)} = 1} \frac{ |g_{f} (i)|}{ \ddz \rz (i)} \rt) \lt( \inf_{i \geq 0} \frac{\ddz \rz(i+1)}{\ddz \rz(i)} \rt) = (m_\rz)^{-1} \sup_{\| f \|_{{\rm Lip} (\rz)} = 1} M_0 (g_f),
\end{align}
}
which yields \Ref{Po0}. \qed

\noindent {\bf Proof of (\ref{Mg0}).} Combining (\ref{I-Q}) and (\ref{I-Q2}), it holds that
\begin{align*}
\sup_{\| f \|_{\Lip} = 1} |g_f (1)| &= \sup_{\| f \|_{\Lip} = 1} - \int_0^\infty e^{-t} \EE \lt[ f(X_t^0 +1) - f(X_t^0) \rt] \d t \\
&= \int_0^\infty e^{-t} \EE \lt[ \ddz \rz (X_t^0) \rt] \d t = (I-Q)^{-1} (\ddz \rz) (0) = - \ddz h_{\rz} (0).
\end{align*}
Hence, using (\ref{stein2}) with $f = \rz$ and $i = 0$, we have
\be\lb{Po0-3}
\sup_{\| f \|_{\Lip} = 1} \frac{|g_f (0)|}{\ddz \rz (0)} = \frac{\pz (\rz) - \rz(0)}{\lz \ddz \rz (0)},
\de
which is (\ref{Mg0}) in Proposition~\ref{remadd}. \qed

\vskip 1 cm

\noindent {\bf Proof of (\ref{Po1}).} Since $\ddz g_f (0) = 0$, we consider $\ddz g_f(i)$ for $i\ge 1$. Using the coupling (\ref{couple_Po}) again, we have
\be\lb{Po1-0}
\ddz g_f (i) = - \int_0^\infty e^{-t} \EE \lt[ \ddz f (X_t^i) - \ddz f(X_t^{i-1}) \rt] \d t = - \int_0^\infty e^{-2t} \EE \lt[ \ddz^2 f (X_t^{i-1}) \rt] \d t, \qd i \geq 1.
\de
This ensures that without loss of generality, we may assume $f(i)=0.$ We now deduce that for any fixed $i\ge 1$, $\sup_{\| f \|_{{\rm Lip} (\rz)} = 1} |\ddz g_f (i)| $ is attained by the function $f_i^* (j) = - |\rz (j) - \rz (i)|$. The argument is exactly the same as in \cite{BX06}, but for the ease of reading, we repeat it here. In fact, \cite[(2.9)]{BX06} says that
 $$\ddz g_f (i)=-\ddz e_{i-1}^+\sum_{j\ge i+1}\pz_jf(j)+\ddz e_i^-\sum_{j\le i-1}\pz_jf(j)+\pz_if(i)(e_{i-1}^++e_{i+1}^-),$$
 and it follows from \Ref{e1} that $\ddz g_f (i)\le \ddz g_{f_i^*} (i)$.

Next, direct computation gives
\be\lb{f2}
\ddz^2 f_i^* (j)=
\left\{
  \begin{array}{ll}
    - \ddz^2 \rz (j), & j \geq i, \\
    \rz(i-1) - \rz (i+1), & j = i-1, \\
    \ddz^2 \rz (j), & j \leq i-2.
  \end{array}
\right.
\de
When $\ddz^2 \rz(i) \geq 0$, $\forall i \geq 1$, we have
\begin{align}\lb{Po1-1}
\sup_{\| f \|_{{\rm Lip} (\rz)} = 1} & | \ddz g_f (i)| = \ddz g_{f_i^*} (i) \notag \\
& = - \int_0^\infty  e^{-2 t} \EE \Big[ - \ddz^2 \rz (X_t^{i-1}) \mathbf{1}_{\{ X_t^{i-1} \geq i \}} + (\rz (X_t^{i-1}) - \rz (X_t^{i-1}+2)) \mathbf{1}_{\{ X_t^{i-1} = i-1 \}} + \ddz^2 \rz (X_t^{i-1}) \mathbf{1}_{\{ X_t^{i-1} \leq i-2 \}} \Big] \d t \notag \\
& = \int_0^\infty  e^{-2 t} \EE \lt[ \ddz^2 \rz (X_t^{i-1})  - 2 (\rz (X_t^{i-1}) - \rz (X_t^{i-1}+1) ) \mathbf{1}_{\{ X_t^{i-1} = i-1 \}} - 2 \ddz^2 \rz (X_t^{i-1}) \mathbf{1}_{\{ X_t^{i-1} \leq i-2 \}} \rt] \d t \notag \\
& \leq \int_0^\infty e^{-2 t} \EE \lt[ \ddz^2 \rz (X_t^{i-1}) \rt] \d t + 2 (\rz (i)- \rz (i-1)) \int_0^\infty e^{-2t} \PP (X_t^{i-1} = i-1) \d t.
\end{align}
It remains to handle the right-hand side of (\ref{Po1-1}).

Firstly, in order to bound $\int_0^\infty e^{-2t} \PP (X_t^{i-1} = i-1) \d t$, we start from the expression (\ref{PoPt}) of the semigroup $P_t$. When $0< \lz \leq 1$, it holds that $(\lz(1 - e^{-t}))^n / (n!) \leq 1$, $\forall n \in \ZZ$, $t \geq 0$. Then by (\ref{PoPt}), we have
\begin{align}\lb{Po1add1}
P_t (i, i) & = e^{-\lz (1- e^{-t})} \sum_{k=0}^{i} \frac{i !}{k! (i-k)!} e^{-k t} (1- e^{-t})^{i-k} \lt( \frac{\lz^{i-k} (1 -  e^{-t})^{i-k}}{(i-k)!} \rt) \notag \\
& \leq e^{-\lz (1- e^{-t})} \sum_{k=0}^{i} \binom{i}{k} e^{-k t} (1- e^{-t})^{i-k} = e^{-\lz (1- e^{-t})}, \qqd t \geq 0.
\end{align}
Hence, we have
\be\lb{Po1-2}
\sup_{i \ge 1} \int_0^\infty e^{-2t} \PP (X_t^{i-1} = i-1) \d t = \int_0^\infty e^{-2t} \PP (X_t^{0} = 0) \d t = \int_0^\infty e^{-2 t} e^{-(\lz - \lz e^{-t})} \d t = \frac{e^{- \lz} + \lz - 1}{\lz^2}.
\de
For $1< \lz < \infty$, \cite[p.~24]{BB92} states that $X_t^{i-1}=X_t^0+Y_t$, where $Y_t\sim {\rm Binomial}(i-1,e^{-t})$ is independent of $X_t^0$ and
\be\lb{PoPt-1}
\PP (X_t^0 = j) = P_t (0, j) = \frac{(\lz(1 - e^{-t}))^j}{j!} e^{- \lz(1 - e^{-t})}, \qqd \forall j \in \ZZ,
\de
hence
\begin{equation}
\PP (X_t^{i-1} = i-1)\le \sup_{j \in \ZZ} \PP (X_t^0 = j),\lb{couplinginequ}
\end{equation}
which ensures
\begin{equation}
\int_0^\infty e^{-2t} \PP (X_t^{i-1} = i-1) \d t \leq \int_0^\infty e^{-2 t} \sup_{j \in \ZZ} \PP (X_t^0 = j) \d t.\label{xia1}
\end{equation}
It is easy to see that (\ref{PoPt-1}) is maximized by the integer-value function $p(t) := \max \{ j \in \ZZ: j \leq \lz - \lz e^{-t} \}$. Obviously, we have $\lt\{ t: p(t) = 0 \rt\} = \lt[ 0 , \log \lz - \log (\lz -1) \rt)$. Applying the following inequality introduced in \cite{XHY97}, which is a more accurate version of Stirling's formula,
$$
r_n \lt( 1+ \frac{1}{12 n} \rt) < n ! < r_n \lt(1+ \frac{1}{12n - 0.5} \rt), \qd n \geq 1, \qd \text{where $r_n := \sqrt{2 \pz n} \lt( \frac{n}{e} \rt)^n$},
$$
then for each $t \geq \log \lz - \log (\lz -1)$, it holds that
\be\lb{Po1-3}
\PP (X_t^0 = p(t)) \leq \frac{1}{\sqrt{2 \pz p(t)} (1+ (12 p(t))^{-1})} \lt( \frac{\lz - \lz e^{-t}}{p(t)} \rt)^{p(t)} e^{p(t) - (\lz - \lz e^{-t})} \leq \frac{1}{\sqrt{2 \pz p(t)} (1+ (12 p(t))^{-1})},
\de
where the last inequality follows from the fact that $( 1+ x / n )^n \leq e^x$, $\forall n \geq 1$, $x \in [0, 1]$. Recall that $\lfloor \lz \rfloor$ is the largest integer less than or equal to $\lz$, for each $1 \leq n \leq \lfloor \lz \rfloor - 1$, we have
$$
\lt\{ t: p(t) = n \rt\} = \lt[ \log \lt( \frac{\lz}{\lz - n} \rt), \log \lt( \frac{\lz}{\lz - n - 1} \rt) \rt), \qd \text{and} \qd \lt\{ t: p(t) = \lfloor \lz \rfloor \rt\} = \lt[ \log \lt( \frac{\lz}{\lz - \lfloor \lz \rfloor} \rt), \infty \rt).
$$
Hence, the integral interval $[0, \infty)$ can be broken down into $\lfloor \lz \rfloor + 1$ parts, and we have
\begin{align}\lb{Po1-4}
\int_0^\infty e^{-2t} \sup_{j \in \ZZ} \PP (X_t^0 = j) \d t & = \int_{\{ t: p(t) = 0 \}} e^{-2 t} e^{- (\lz - \lz e^{-t})} \d t + \sum_{n=1}^{\lfloor \lz \rfloor -1} \frac{1}{\sqrt{2 \pz n} (1+ (12n)^{-1})} \int_{\{ t : p(t) =n \}} e^{-2t} \d t \notag \\
& \qd \qd + \frac{1}{\sqrt{2 \pz \lfloor \lz \rfloor} (1+ (12\lfloor \lz \rfloor)^{-1})} \int_{ \{ t: p(t) = \lfloor \lz \rfloor \}} e^{-2t} \d t \notag \\
&= \frac{(e - 1)(\lz-1) +1 }{\lz^2 e} + \sum_{n=1}^{\lfloor \lz \rfloor -1} \lt( \frac{ 12 \sqrt{n} (\lz -n) - 6 \sqrt{n} }{ \sqrt{2 \pz} \lz^2 (12 n + 1)} \rt) + \frac{6 \sqrt{ \lfloor \lz \rfloor} (\lz - \lfloor \lz \rfloor)^2 }{\sqrt{2 \pz} \lz^2 (12 \lfloor \lz \rfloor + 1)}.
\end{align}
Secondly, for the estimate of $\int_0^\infty e^{-2 t} \EE \lt[ \ddz^2 \rz (X_t^{i-1}) \rt] \d t$, we use the coupling (\ref{couple_Po}) and the formulae~(\ref{I-Q}), {(\ref{I-Q2})}
to obtain
\begin{align}\lb{Po1-5xia1}
\int_0^\infty e^{-2t} \EE \lt[ \ddz^2 \rz(X_t^{i-1}) \rt] \d t & = \int_0^\infty e^{-t} \EE \lt[ \ddz \rz (X_t^i) - \ddz \rz (X_t^{i-1}) \rt] \d t \notag \\
&= \int_0^\infty e^{-t} \lt[ P_t (\ddz \rz) (i) - P_t (\ddz \rz) (i-1) \rt] \d t \notag \\
&= (I-Q)^{-1} (\ddz \rz) (i) - (I-Q)^{-1} (\ddz \rz) (i-1)\notag\\
&= - \ddz^2 h_{\rz} (i-1).
\end{align}

Combining \Ref{Po1-1}, \Ref{Po1-2}, \Ref{xia1}, \Ref{Po1-4} and \Ref{Po1-5xia1}, we have
\begin{align}\lb{Po1-6}
\sup_{\| f \|_{{\rm Lip} (\rz)} = 1}M_1 (g_f) & = \sup_{i \geq 1} \frac{\ddz g_{f_i^*} (i)}{\ddz\rz (i)} \leq \sup_{i \geq 1} \frac{|\ddz h_{\rz} (i) - \ddz h_{\rz} (i-1)|}{ \ddz \rz (i) } + 2 m_\rz\sup_{i \geq 1} \int_0^\infty e^{-2t} \PP (X_t^{i-1} = i-1) \d t \notag \\
& \leq m_\rz\lt\| \ddz h_{\rz} \rt\|_{\Lip} + 2m_\rz \Xi_1 (\lz),
\end{align}
where $\Xi_1(\lz)$ is defined in (\ref{lz1}). \qed

\vskip1cm

\begin{rem}If $\rz(i) = \rz_1 (i) = i$ and $0< \lz \leq 1$, the estimate of $\sup_{\| f \|_{{\rm Lip} (\rz)} = 1} M_1 (g_f)$ is sharp.
\end{rem}

In fact, since $\ddz \rz_1 (i) =1$ and $\ddz^2 \rz_1 (i) = 0$, $m_\rz=1$, using (\ref{Po1-1}) and (\ref{Po1-2}), we have
\be\lb{Po1add2}
\sup_{\| f \|_{{\rm Lip} (\rz)} = 1} M_1 (g_f) = 2 \sup_{i \geq 1} \int_0^\infty e^{-2t} \PP \lt( X_t^{i-1} = i-1\rt) \d t = 2 \int_0^\infty e^{-2t} \PP \lt( X_t^{0} = 0 \rt) \d t = \frac{2 (e^{-\lz} + \lz -1) }{\lz^2}.
\de

\vskip1cm

\noindent {\bf Proof of (\ref{XiaPo1}).} When $\ddz^2 \rz (i) \leq 0$, $\forall i \in \ZZ$, one can repeat the proof of \Ref{Po1} but replace \Ref{Po1-1} with
\begin{align}\lb{XiaPo1-1}
\sup_{\| f \|_{\Lip} = 1} & | \ddz g_f (i)| = \ddz g_{f_i^*} (i) \notag \\
& = \int_0^\infty  e^{-2 t} \EE \Big[  \ddz^2 \rz (X_t^{i-1}) \mathbf{1}_{\{ X_t^{i-1} \geq i \}} - (\rz (X_t^{i-1}) - \rz (X_t^{i-1}+2)) \mathbf{1}_{\{ X_t^{i-1} = i-1 \}} - \ddz^2 \rz (X_t^{i-1}) \mathbf{1}_{\{ X_t^{i-1} \leq i-2 \}} \Big] \d t \notag \\
& = \int_0^\infty  e^{-2 t} \EE \lt[ -\ddz^2 \rz (X_t^{i-1})  + (\ddz^2 \rz (X_t^{i-1}) + \rz (X_t^{i-1}+2) - \rz (X_t^{i-1}) ) \mathbf{1}_{\{ X_t^{i-1} = i-1 \}} + 2 \ddz^2 \rz (X_t^{i-1}) \mathbf{1}_{\{ X_t^{i-1} \geq i \}} \rt] \d t \notag \\
& \leq \int_0^\infty e^{-2 t} \EE \lt[ -\ddz^2 \rz (X_t^{i-1}) \rt] \d t + 2 \ddz\rz (i) \int_0^\infty e^{-2t} \PP (X_t^{i-1} = i-1) \d t, \mbox{ for }i\ge 1,
\end{align}
and then
\begin{align*}
\sup_{\| f \|_{{\rm Lip} (\rz)} = 1}M_1 (g_f) & = \sup_{i \geq 1} \frac{\ddz g_{f_i^*} (i)}{\ddz\rz (i)} \leq \sup_{i \geq 1} \frac{|\ddz h_{\rz} (i) - \ddz h_{\rz} (i-1)|}{ \ddz \rz (i) } + 2 \sup_{i \geq 1} \int_0^\infty e^{-2t} \PP (X_t^{i-1} = i-1) \d t \notag \\
& \leq m_\rz\lt\| \ddz h_{\rz} \rt\|_{\Lip} + 2 \Xi_1 (\lz). \qed
\end{align*}

\vskip 1 cm

\noindent {\bf Proof of (\ref{Mg2}).} Since $\ddz g_f(0) = 0$, we have $\ddz^2 g_f(0)=\ddz g_f(1)$. Using \Ref{Po1-1},  \Ref{XiaPo1-1} and \Ref{Po1-5xia1} with $i=1$, we obtain
\be\lb{XiaPo2-01}
\sup_{\| f\|_{\Lip} = 1} \frac{|\ddz^2 g_f (0)|}{\ddz \rz (0)} = \frac{|\ddz^2 h_{\rz} (0)|}{\ddz \rz (0)} +
\left\{
  \begin{array}{ll}
    \dps\frac{2(e^{-\lz} + \lz -1)}{\lz^2}, & \hbox{when $\ddz^2 \rz (\cdot) \geq 0$;} \\
    \dps\frac{2 \ddz \rz (1) (e^{-\lz} + \lz -1)}{\ddz \rz (0) \lz^2}, & \hbox{when $\ddz^2 \rz (\cdot) \leq 0$.}
  \end{array}
\right.
\de
It follows from (\ref{stein2}) with $f = \rz$ and $i = 0, 1$ that
\be\lb{Po2-add}
\frac{|\ddz^2 h_{\rz} (0)|}{\ddz \rz (0)} = \lt| \frac{\rz(1) - \pz (\rz) + \ddz h_\rz (0)}{\lz \ddz \rz (0)} - \frac{\ddz h_\rz (0)}{\ddz \rz (0)} \rt| = \lt| \frac{1}{\lz} + \frac{\rz(0) - \pz (\rz)}{\lz^2 \ddz \rz (0)} \rt|.
\de
Hence, (\ref{Mg2}) in Proposition~\ref{remadd} is implied by \Ref{XiaPo2-01} and \Ref{Po2-add}. \qed

\vskip 1cm

\noindent {\bf Proof of (\ref{Po2}).} Now, we can focus on $\ddz^2 g_f(i)$ for $i \geq 1$. Combining (\ref{couple_Po}) and (\ref{Po1-0}), we have
\be\lb{XiaPo2-1}\ddz^2 g_f (i) = - \int_0^\infty e^{-2t} \EE \lt[ \ddz^2 f (X_t^i) - \ddz^2 f (X_t^{i-1}) \rt] \d t = - \int_0^\infty e^{-3 t} \EE \lt[ \ddz^3 f (X_t^{i-1}) \rt] \d t.\de
Hence, without loss of generality, we may again take $f(i)=0$.
As in \cite{BX06}, we argue that $\sup_{\| f \|_{\Lip} = 1} |\ddz^2 g_f (i)|$ is achieved by the function $f_i^{\triangle}$ defined as
$$
f_i^{\triangle} (j) =
\left\{
  \begin{array}{ll}
    \rz(i) - \rz (j), & 1 \leq j \leq i, \\
    2 \rz(i+1) - \rz(i) - \rz(j), & j \geq i+1.
  \end{array}
\right.
$$
For the sake of completeness, we recall the proof of \cite{BX06} here. In fact, \cite[(2.18)]{BX06} states
$$\ddz^2 g_f (i)=-\ddz^2 e_{i-1}^+\sum_{j\ge i+2}(f(j)-f(i+1))\pz_j+\ddz^2 e_i^-\sum_{j\le i-1}\pz_jf(j)+f(i+1)r_i.$$
Hence, we can see from \Ref{e2} and \Ref{xiari} that $\ddz^2 g_f (i)\le \ddz^2 g_{f_i^{\triangle}} (i)$. This, together with \Ref{XiaPo2-1}, ensures
\be\lb{Po2-1}
\sup_{\| f \|_{\Lip} = 1} \ddz^2 g_f (i) = \ddz^2 g_{f_i^{\triangle}} (i) = - \int_0^\infty e^{-3 t} \EE \lt[ \ddz^3 f_i^{\triangle} (X_t^{i-1}) \rt] \d t,
\de
thus, it suffices to estimate $\EE \lt[ \ddz^3 f_i^{\triangle} (X_t^{i-1}) \rt]$. Since
$$
\ddz^3 f_i^{\triangle} (j) =
\left\{
  \begin{array}{ll}
    - \ddz^3 \rz (j), & \hbox{$j \leq i-3$ or $j \geq i+1$,} \\
      \rz(i+1) + \rz (i) - 3 \rz (i-1) + \rz (i-2), & j= i-2, \\
    -\rz (i+2) - \rz(i+1) + \rz (i) + \rz (i-1), & j= i-1, \\
    - \rz (i+3) + 3 \rz(i+2) - \rz(i+1) - \rz(i), & j=i, \\
  \end{array}
\right.
$$
we obtain
\begin{align}\lb{delta2xiaadd1}
\EE \lt[ \ddz^3 f_i^{\triangle} (X_t^{i-1}) \rt] &= - \sum_{j \leq i-3} \ddz^3 \rz(j) \PP (X_t^{i-1} = j) - \sum_{j \geq i+1} \ddz^3 \rz(j) \PP (X_t^{i-1} = j) \notag\\
& \qqd + \lt[ \rz(i+1) + \rz (i) - 3 \rz (i-1) + \rz (i-2) \rt] \PP (X_t^{i-1} = i-2) \notag\\
& \qqd + \lt[ -\rz (i+2) - \rz(i+1) + \rz (i) + \rz (i-1) \rt] \PP (X_t^{i-1} = i-1) \notag\\
& \qqd + \lt[ - \rz (i+3) + 3 \rz(i+2) - \rz(i+1) - \rz(i) \rt] \PP (X_t^{i-1} = i) \notag\\
&= - \EE \lt[ \ddz^3 \rz (X_t^{i-1}) \rt] + 2 \ddz \rz(i) \lt[ \PP (X_t^{i-1} = i-2) - 2 \PP (X_t^{i-1} = i-1) + \PP (X_t^{i-1} = i) \rt].
\end{align}
Combining \Ref{Po2-1} and \Ref{delta2xiaadd1} gives
\be\lb{Po2-2}
\sup_{\| f \|_{\Lip} = 1} M_2 (g_f) \leq \sup_{i \geq 1} \frac{\int_0^\infty e^{-3t} \EE \lt[ \ddz^3 \rz (X_t^{i-1}) \rt] \d t }{\ddz \rz (i)} + 4 \int_0^\infty e^{-3t} \PP (X_t^{i-1} = i-1) \d t.
\de
For the first item of (\ref{Po2-2}), by (\ref{couple_Po}) and (\ref{Po1-5xia1}), we have
\begin{align}\lb{Po2-3}
\sup_{i \geq 1} & \lt( \ddz \rz(i) \rt)^{-1} \lt[ \int_0^\infty e^{-2t} \EE \lt[ \ddz^2 \rz (X_t^{i}) \rt] \d t - \int_0^\infty e^{-2t} \EE \lt[ \ddz^2 \rz (X_t^{i-1}) \rt] \d t \rt] \notag \\
& \qd = \sup_{i \geq 1} \frac{\lt| \ddz^2 h_{\rz} (i) - \ddz^2 h_{\rz} (i-1) \rt|}{\ddz \rz (i)} \notag \\
& \qd \le m_\rz\lt\| \ddz^2 h_{\rz} \rt\|_{\Lip}.
\end{align}
For the second item of (\ref{Po2-2}), using the estimate given in (\ref{couplinginequ}), we have
\be\lb{Po2-4}
4 \int_0^\infty e^{-3 t} \lt[ \PP (X_t^{i-1} = i-1) \rt] \d t \leq 4 \int_0^\infty e^{-3 t} \lt( \sup_{j \in \ZZ} P_t (0, j) \rt) \d t.
\de
To bound $\int_0^\infty e^{-3 t} \lt( \sup_{j \in \ZZ} P_t (0, j) \rt) \d t$, we use the same argument as that in the proof of (\ref{Po1}). When $0< \lz \leq 1$, we have
\be\lb{Po2-5}
\int_0^\infty e^{-3 t} \lt( \sup_{j \in \ZZ} P_t (0, j) \rt) \d t \leq \int_0^\infty e^{-3 t} e^{- (\lz - \lz e^{-t})} \d t = \frac{(\lz -1)^2 -2 e^{-\lz} +1}{\lz^3}.
\de
When $1< \lz < \infty$, using the same notation $p(t)$ introduced in the proof of (\ref{Po1}), we have
\begin{align}\lb{Po2-6}
\int_0^\infty e^{-3 t} \lt( \sup_{j \in \ZZ} P_t (0, j) \rt) \d t & = \int_{\{t: p(t) =0 \}} e^{-3 t} e^{- (\lz - \lz e^{-t})} \d t + \sum_{n=1}^{\lfloor \lz \rfloor -1} \frac{1}{\sqrt{2 \pz n} (1+ (12n)^{-1})} \int_{ \{t: p(t) = n \}} e^{-3t} \d t \notag \\
& \qd + \frac{1}{\sqrt{2 \pz \lfloor \lz \rfloor} (1+ (12\lfloor \lz \rfloor)^{-1})} \int_{\{t: p(t) = \lfloor \lz \rfloor \}} e^{-3t} \d t \notag \\
&= \frac{\lz^2 (e - 1) -2 \lz (e-2) + 2e -5}{\lz^3 e} + \sum_{n=1}^{\lfloor \lz \rfloor -1} \lt( \frac{4 \sqrt{n} (3 (\lz - n)^2 - 3(\lz - n) + 1)}{ \sqrt{2 \pz} \lz^3 (12n + 1)} \rt) \notag \\
& \qd \qd + \frac{4 \sqrt{\lfloor \lz \rfloor} ( \lz - \lfloor \lz \rfloor)^3}{\sqrt{2 \pz} \lz^3 (12 \lfloor \lz \rfloor + 1)}.
\end{align}
Hence, by (\ref{Po2-3}) -- (\ref{Po2-6}), we have
\be\lb{Po2-7}
\sup_{\| f \|_{\Lip} = 1} M_2 (g_f) \leq m_\rz\lt\| \ddz^2 h_{\rz} \rt\|_{\Lip} + 4  \Xi_2 (\lz),
\de
where $\Xi_2 (\lz)$ is defined in (\ref{lz2}).

Finally, we use another method to bound $\ddz^2 g_{f_i^{\triangle}} (i)$, which is different from (\ref{Po2-7}). Note that by the representation of $g_f$ in (\ref{Po0-0}), we have from (\ref{e0}) that
$$
g_f (i) = e_i^- \sum_{j=0}^{i-1} \pz_j f(j) - e_{i-1}^+ \sum_{j=i}^\infty \pz_j f(j),
$$
which means that $g_f$ has linear property with respect to $f$. Moreover,
\begin{align}\lb{Po2-8}
\ddz^2 g_f (i) & = g_f (i+2) - 2 g_f (i+1) + g_f (i) \notag \\
&= \lt( \ddz^2 e_i^- \rt) \sum_{j=0}^{i-1} \pz_j f(j) - \lt( \ddz^2 e_{i-1}^+ \rt) \sum_{j = i+2}^\infty \pz_j f (j) \notag \\
& \qd + \lt( 2e_{i}^+ - e_{i-1}^+ + e_{i+2}^- \rt) \pz_{i+1} f (i+1) + \lt( e_{i+2}^- - 2e_{i+1}^- - e_{i-1}^+ \rt) \pz_{i} f (i).
\end{align}
Given any $i \ge 1$, define $\fz_i (j) = \rz (i)- \rz (j)$, for $ j \in \ZZ$, it follows from \Ref{XiaPo2-1} that
\be\lb{Po2-9}
\ddz^2 g_{\fz_i} (i) = \int_0^\infty e^{-3 t} \EE \lt[ \ddz^3 \rz (X_t^{i-1}) \rt] \d t,
\de
and
$$
f_i^{\triangle} (j) - \fz_i (j) =
\left\{
  \begin{array}{ll}
    0, & 1 \leq j \leq i, \\
    2 \ddz\rz (i), & i+1 \leq j < \infty.
  \end{array}
\right.
$$
Using (\ref{Po2-8}) directly, we have
\begin{align}\lb{Po2-10}
\ddz^2 g_{f_i^{\triangle} - \fz_i} (i) & = -2 \lt( \ddz^2 e_{i-1}^+ \rt) \ddz\rz (i) \RF (i+2) + 2 \ddz\rz (i) \lt( 2e_i^+ - e_{i-1}^+ + e_{i+2}^- \rt) \pz_{i+1} \notag \\
& = -2 \lt( \ddz^2 e_{i-1}^+ \rt) \ddz\rz (i) \RF (i+1) + 2 \ddz\rz (i) \lt( e_{i+1}^+ + e_{i+2}^- \rt) \pz_{i+1} \notag \\
& \leq 2 \ddz\rz (i)/\lz,
\end{align}
where the last inequality is due to (\ref{e2}) and $\pz_{i+1} \lt( e_{i+1}^+ + e_{i+2}^- \rt) = \lz^{-1}$. By \Ref{Po1-5xia1},  (\ref{Po2-9}), (\ref{Po2-10}) and the linear property of $g_f$, we obtain
\begin{align}\lb{Po2-11}
\sup_{\| f \|_{\Lip} = 1} M_2 (g_f) & \leq \sup_{i \geq 1} \frac{ \ddz^2 g_{f_i^{\triangle}} (i)}{\ddz \rz (i)}  \le \sup_{i \geq 1} \frac { \ddz^2 g_{\fz_i} (i)}{\ddz \rz (i)} + \sup_{i \geq 1} \frac{\ddz^2 g_{f_i^{\triangle} - \fz_i} (i)}{\ddz \rz (i)} \notag \\
&\leq m_\rz\lt\| \ddz^2 h_{\rz} \rt\|_{\Lip} + \frac{2}{\lz} .
\end{align}
Combining (\ref{Po2-7}) and (\ref{Po2-11}), we obtain
$$
\sup_{\| f \|_{\Lip} = 1} M_2 (g_f) \leq m_\rz\lt\| \ddz^2 h_{\rz} \rt\|_{\Lip} + 2 [(2 \Xi_2 (\lz)) \wedge \lz^{-1}],
$$
and the proof of Theorem \ref{thmPo} is complete. \deprf

\vskip 1 cm

\noindent {\bf Proof of Proposition \ref{prop}.} (1). Let $\rz_p (i) = i^p$, $p \geq 1$. Obviously, for each $i \in \ZZ$, it holds that $\ddz \rz_p (i) \geq 0$, $\ddz^2 \rz_p (i) \geq 0$ and
\be\lb{ex1-1}
\pz (\rz_p) = \lz \sum_{i \geq 0} \pz_i (i+1)^{p-1} = \lz \sum_{i \geq 0} \pz_i \sum_{k=0}^{p-1} \binom{p-1}{k} i^k = \lz \sum_{k=0}^{p-1} \binom{p-1}{k} \pz (\rz_k).
\de
Note that $h_p$ is the solution to the Stein equation (\ref{Poisson}), that means $h_p (i) = Q^{-1} (\rz_p - \pz (\rz_p)) (i)$, $\forall i \geq 0$. When $p =1$ and $i \geq 1$, it holds that $Q \rz_1 (i) = \lz - i $, which implies that $h_1 (i) = Q^{-1} (\rz_1 - \pz (\rz_1)) (i) = -i$, $i \geq 1$. For $i =0$, since $Q h_1 (0) = \lz (h_1 (1) - h_1 (0)) = - \lz$, we have $h_1 (0) = 0$. When $p \geq 2$ and $i \geq 1$, we have
\begin{align}
Q (i^p) & = \lz \lt((i+1)^p - i^p\rt) + i \lt( (i-1)^p - i^p\rt) = - p i^p + \lz + \sum_{k=1}^{p-1} \binom{p}{k} i^k \lt[ \lz + \frac{k (-1)^{p-k+1}}{p-k+1} \rt] \nonumber\\
&= -p \lt(i^p - \pz (\rz_p) \rt) + \sum_{k=1}^{p-1} \binom{p}{k} \lt(i^k - \pz (\rz_k) \rt) \lt[ \lz + \frac{k (-1)^{p-k+1}}{p-k+1} \rt],\label{prop1.7xia1}
\end{align}
where the last equality is based on the following observation: with $\eta\sim \pz$, $j=k-1$, using \Ref{Steinop}, we have
\begin{align*}
&  - p\pz (\rz_p) +\lz+ \sum_{k=1}^{p-1} \binom{p}{k}  \pz (\rz_k) \lt[ \lz + \frac{k (-1)^{p-k+1}}{p-k+1} \rt] \\
& = -p\EE(\eta^p)+\sum_{j=0}^{p-2}(-1)^{p-j}{p\choose j}\EE(\eta^{j+1})+\lz\sum_{k=0}^{p-1}{p\choose k}\EE(\eta^k)\\
& = -p\EE(\eta^p)+\EE\left\{\eta\left((\eta-1)^p-\eta^p+p\eta^{p-1}\right)\right\}+\lz\EE\left((\eta+1)^p-\eta^p\right) \\
&= -p\EE(\eta^p)+\lz\EE(\eta^p)-\lz\EE\left((\eta+1)^p\right)+p\EE(\eta^{p})+\lz\EE\left((\eta+1)^p-\eta^p\right) \\
& = 0.
\end{align*}
Hence, applying $Q^{-1}$ to both sides of \Ref{prop1.7xia1}, by the definition of $h_p (i)$, we obtain
$$
h_p (i) = - \frac{1}{p} i^p + \frac{1}{p} \sum_{k=1}^{p-1} \binom{p}{k} h_k (i) \lt[ \lz + \frac{k (-1)^{p-k+1}}{p-k+1} \rt], \qqd i \geq 1.
$$
Similarly, since $Q h_p (0) = \lz (h_p (1) - h_p (0)) = - \pz (\rz_p)$, we have $h_p (0) = h_p (1) + \lz^{-1} \pz(\rz_p)$.

In particular, when $p=2$, we have $m_\rz=1$, $\pz (\rz_2) = \lz^2 + \lz$ and
$$
h_2 (i) = - \frac{1}{2} i^2 - \frac{1}{2} (2 \lz + 1) i, \qqd i \geq 1.
$$
According to the expression of $\| (-Q)^{-1} \|_{\Lip}$ in (\ref{Po0-2}), we have
$$
\lt\| Q^{-1} \rt\|_{\text{\rm Lip}(\rz_2)} = \sup_{i \geq 1} \frac{|h_2 (i) - h_2 (i-1)|}{\rz_2 (i) - \rz_2 (i-1)} = \frac{1}{2} \sup_{i \geq 1} \lt( 1+ \frac{2 \lz +1}{2 i -1} \rt) = 1+ \lz.
$$
Since
$$
\ddz h_2 (i) = -i -\lz -1, \ \ \ddz^2 h_2 (i) = -1 \ \ \text{and} \ \ \ddz^3 h_2 (i) =0, \qqd i \geq 0,
$$
we have
$$
\lt\| \ddz h_2 \rt\|_{\text{\rm Lip}(\rz_2)} = \sup_{i \geq 1} \frac{|\ddz^2 h_2 (i-1)|}{\ddz \rz (i-1)} = \sup_{i \geq 1} \frac{1}{2i -1} = 1 \qd \text{and} \qd \lt\| \ddz^2 h_2 \rt\|_{\text{\rm Lip}(\rz_2)} = \sup_{i \geq 1} \frac{|\ddz^3 h_2 (i-1)|}{\ddz \rz (i-1)} = 0.
$$
Finally, according to Theorem \ref{thmPo}, we obtain the estimate (\ref{rz2}).

(2). Let
$$
\rz_{1/2} (i) = \lz + \sqrt{i} - \frac{\lz}{\sqrt{i+1}} \geq 0, \qqd \forall i \in \ZZ.
$$
Then we have
$$
\pz (\rz_{1/2}) = \lz + \sum_{i \geq 1} \frac{\lz}{\sqrt{i}} e^{-\lz} \frac{\lz^{i-1}}{(i-1)!} - \sum_{i \geq 0} \frac{\lz}{\sqrt{i+1}} e^{-\lz} \frac{\lz^{i}}{i!} = \lz.
$$
For each $i \in \ZZ$,
\begin{align*}
\ddz \rz_{1/2} (i) &= \frac{1}{\sqrt{i} + \sqrt{i+1}} + \lz \lt( \frac{\sqrt{i+2} - \sqrt{i+1}}{\sqrt{(i+1)(i+2)}}\rt) > 0; \\
\ddz^2 \rz_{1/2} (i) &= \frac{(-1) \lt( \sqrt{i+2} - \sqrt{i} \rt)}{\lt( \sqrt{i+1} + \sqrt{i} \rt) \lt(\sqrt{i+1} + \sqrt{i+2} \rt)} - \lz \lt( \ddz^2 (i+1)^{- \frac{1}{2}} \rt) < 0.
\end{align*}
Moreover, it is easy to demonstrate that
$$
h_{1/2} (0) = 0, \qd h_{1/2} (i) = - \sum_{k=1}^i \frac{1}{\sqrt{k}}, \qqd i \geq 1,
$$
satisfies the Stein equation (\ref{Poisson}), i.e. $Q h_{1/2} = \rz_{1/2} - \pz (\rz_{1/2})$, and
$$
\ddz h_{1/2} (i) = \frac{-1}{\sqrt{i+1}}, \qd \ddz^2 h_{1/2} (i) = \frac{1}{(i+1) \sqrt{i+2}+ (i+2) \sqrt{i+1}}, \qqd i \geq 1.
$$
Here, we introduce an auxiliary function $\fz (i)$,
\be\lb{ex2-0}
\fz (i) = \frac{\sqrt{i+1} + \sqrt{i}}{\sqrt{i} + \sqrt{i-1}} = \frac{\sqrt{1+ i^{-1}} + 1}{1+ \sqrt{1- i^{-1}}}, \qqd i \geq 1.
\de
It is easy to verify that $\fz(i) \geq 1$ and $\fz(i)$ is decreasing for each $i \geq 1$.

Firstly, we consider $\| (-Q)^{-1} \|_{\text{\rm Lip}(\rz_{1/2})}$. By its definition, we have
\begin{align}\lb{ex2-1}
\| (-Q)^{-1} \|_{\text{\rm Lip}(\rz_{1/2})} & = \sup_{i \geq 1} \frac{| h_{1/2} (i) - h_{1/2} (i-1) |}{\rz_{1/2} (i) - \rz_{1/2} (i-1)} = \sup_{i \geq 1} \frac{1}{i - \frac{\lz \sqrt{i}}{\sqrt{i+1}} - \sqrt{i (i-1)} + \lz} \notag \\
&= \sup_{i \geq 1} \frac{1}{\lz \lt(1- \sqrt{1- (i+1)^{-1}} \rt) + \sqrt{i} \fz(i) / \lt(\sqrt{i+1} + \sqrt{i} \rt)}.
\end{align}
Since
\begin{align*}
& 1- \sqrt{1- (i+1)^{-1}} \qd \text{is decreasing and approaching to $0$ as $i \to \infty$}, \\
& \frac{\sqrt{i} \fz(i)}{\sqrt{i+1} + \sqrt{i}} = \frac{1}{1+\sqrt{1-i^{-1}}} \qd \text{is decreasing and approaching to $1/2$ as $i \to \infty$},
\end{align*}
the maximum of (\ref{ex2-1}) is attained at $i \to \infty$. Hence,
$$
\sup_{\| f \|_{{\rm Lip} (\rz_{1/2})} = 1} M_0 (g_f) { \leq m_{\rz}} \| (-Q)^{-1} \|_{\text{\rm Lip}(\rz_{1/2})} = 2 {m_{\rz}}.
$$
{To calculate $m_{\rz} = \sup_{i \geq 0} \frac{\ddz \rz_{1/2}(i)}{\ddz \rz_{1/2} (i+1)}$, we define
$$
F(i) := \frac{\ddz \rz_{1/2}(i)}{\ddz \rz_{1/2} (i+1)} = \frac{\fz(i+2) \sqrt{i+3}}{\sqrt{i+1}} \lt( \frac{\lz + \fz(i+1) \sqrt{(i+1)(i+2)}}{\lz + \fz (i+2) \sqrt{(i+2)(i+3)}} \rt), \qd i \geq 0.
$$
Using the ratio formula, we have
$$
F(i) \geq \frac{\fz(i+2) \sqrt{i+3}}{\sqrt{i+1}} \lt(1 \wedge \frac{\fz(i+1) \sqrt{i+1}}{\fz(i+2) \sqrt{i+3}} \rt) = \lt( \fz(i+2) \sqrt{1 + 2/(i+1)} \rt) \wedge \fz(i+1), \qd i \geq 0.
$$
Note that $\fz(i+2) \sqrt{1 + 2/(i+1)}$ and $\fz(i+1)$ are decreasing for each $i \geq 0$, which implies that $m_{\rz} = \sup_{i \geq 0} F(i) \geq \lt( \fz(2) \sqrt{3}\rt) \wedge \fz(1) = \fz(2) \sqrt{3}$. Using the ratio formula again, for each $i \geq 1$, we have
$$
\sup_{i \geq 1} F(i) \leq \lt( \sup_{i \geq 1} \fz(i+2) \sqrt{1 + 2/(i+1)} \rt) \vee \lt( \sup_{i \geq 1} \fz(i+1) \rt) = \fz(3) \sqrt{2} < \fz(2) \sqrt{3}.
$$
Hence,
$$
m_{\rz} = F(0) \vee \lt( \sup_{i \geq 1} F(i) \rt) = F(0) \vee \lt( \sqrt{3} \fz(2) \rt) = \frac{\sqrt{3} \lt( \sqrt{2} + \sqrt{2}\lz - \lz \rt)}{\sqrt{3} \lt(2-\sqrt{2} \rt) + \lz \lt( \sqrt{3}-\sqrt{2} \rt)}.
$$
}

Secondly, we consider $\| \ddz h_{1/2} \|_{\text{\rm Lip}(\rz_{1/2})}$. Supplement the value of $\ddz h_{1/2} (i)$ at $i =0$ by $\ddz h_{1/2} (0) = h_{1/2} (1) - h_{1/2} (0) = -1$. Again, we begin with the definition
\begin{align}\lb{ex2-2}
\|\ddz h_{1/2} \|_{\text{\rm Lip}(\rz_{1/2})} &= \sup_{i \geq 1} \frac{|\ddz h_{1/2} (i) - \ddz h_{1/2} (i-1)|}{\rz_{1/2} (i) - \rz_{1/2} (i-1)} \notag \\
&= \lt( \sup_{i \geq 2} \frac{1/ \lt( i \sqrt{i+1} + (i+1) \sqrt{i} \rt)}{\sqrt{i} - \lz/ \sqrt{i+1} - \sqrt{i-1} + \lz/ \sqrt{i}} \rt) \vee \lt( \frac{|\ddz h_{1/2} (1) - \ddz h_{1/2} (0)|}{\rz_{1/2} (1) - \rz_{1/2} (0)}\rt) \notag \\
& = \lt( \sup_{i \geq 2} \frac{1}{\lz + \sqrt{i(i+1)} \fz (i)} \rt) \vee \lt( \frac{1}{\lz + 2 + \sqrt{2}}\rt).
\end{align}
Note that $\fz (i) \geq 1$ for each $i \geq 1$, then we have
\be\lb{ex2-3}
\lz + \sqrt{6} \fz(2) \leq \lz + 2 + \sqrt{2} \leq \lz + \sqrt{i(i+1)} \leq \lz + \sqrt{i(i+1)} \fz(i), \qqd \forall i \geq 3.
\de
Hence,
\begin{align}\lb{ex2-4}
\|\ddz h_{1/2} \|_{\text{\rm Lip}(\rz_{1/2})} & = \lt( \sup_{i \geq 3} \frac{1}{\lz + \sqrt{i (i+1)} \fz(i)} \rt) \vee \lt( \frac{1}{\lz + \sqrt{6} \fz(2)} \rt) \vee \lt( \frac{1}{\lz + 2 + \sqrt{2}} \rt) \notag \\
& = \frac{1}{\lz + \sqrt{6} \fz(2)} = \frac{1}{\lz + (\sqrt{2} + \sqrt{3}) (2 \sqrt{3} - \sqrt{6})}.
\end{align}
According to Theorem \ref{thmPo}, we obtain
$$
\sup_{\| f \|_{{\rm Lip} (\rz_{1/2})} = 1} M_1 (g_f) \leq \frac{{m_{\rz}}}{\lz + (\sqrt{2} + \sqrt{3}) (2 \sqrt{3} - \sqrt{6})} + 2 \Xi_1 (\lz).
$$

Finally, we consider $\| \ddz^2 h_{1/2} \|_{\text{\rm Lip}(\rz_{1/2})}$. Similarly, we supplement the value of $\ddz^2 h_{1/2} (i)$ at $i =0$ as $\ddz^2 h_{1/2} (0) = \ddz h_{1/2} (1) - \ddz h_{1/2} (0) = (\sqrt{2}-1)/ \sqrt{2}$. By definition,
\begin{align*}
\| \ddz^2 h_{1/2} \|_{\text{\rm Lip}(\rz_{1/2})} & = \lt( \sup_{i \geq 2} \frac{\ddz \lt[ 1/ \lt(i \sqrt{i+1} + (i+1) \sqrt{i} \rt) \rt]}{\sqrt{i} - \lz/ \sqrt{i+1} - \sqrt{i-1} + \lz/ \sqrt{i}} \rt) \vee \lt( \frac{|\ddz^2 h_{1/2} (1) - \ddz^2 h_{1/2} (0)|}{\rz_{1/2} (1) - \rz_{1/2} (0)} \rt) \\
&= \lt( \sup_{i \geq 2} \frac{1 - \frac{\sqrt{i}}{\sqrt{i+2} \fz(i+1)}}{\lz + \fz(i) \sqrt{i(i+1)}} \rt) \vee \lt( \frac{\lt(\sqrt{2} + 1\rt) \lt(2+ \sqrt{2} - \sqrt{6}/3\rt)}{\lz +2 + \sqrt{2}} \rt).
\end{align*}
Since {$\frac{\sqrt{i}}{\sqrt{i+2} \fz (i+1)}$} is increasing, using (\ref{ex2-3}) again, we have
$$
\sup_{i \geq 3} \frac{1 - \frac{\sqrt{i}}{\sqrt{i+2} \fz(i+1)}}{\lz + \fz(i) \sqrt{i(i+1)}} \leq \sup_{i \geq 3} \frac{1 - \frac{\sqrt{i}}{\sqrt{i+2} \fz(i+1)}}{\lz + \sqrt{6} \fz (2)} = \frac{1 - \frac{\sqrt{3}}{\sqrt{5} \fz(4)}}{\lz + \sqrt{6} \fz (2)} \leq \frac{1 - \frac{\sqrt{2}}{2 \fz(3)}}{\lz + \sqrt{6} \fz (2) }.
$$
Hence,
\begin{align*}
\| \ddz^2 h_{1/2} \|_{\text{\rm Lip}(\rz_{1/2})} & = \lt( \sup_{i \geq 3} \frac{1 - \frac{\sqrt{i}}{\sqrt{i+2} \fz(i+1)}}{\lz + \fz(i) \sqrt{i(i+1)}} \rt) \vee \lt( \frac{1 - \frac{\sqrt{2}}{2 \fz(3)}}{\lz + \sqrt{6} \fz (2) } \rt) \vee \lt( \frac{\lt(\sqrt{2} + 1\rt) \lt(2+ \sqrt{2} - \sqrt{6}/3\rt)}{\lz +2 + \sqrt{2}} \rt) \\
&= \lt( \frac{1 - \frac{\sqrt{2}}{2 \fz(3)}}{\lz + \sqrt{6} \fz (2) } \rt) \vee \lt( \frac{\lt(\sqrt{2} + 1\rt) \lt(2+ \sqrt{2} - \sqrt{6}/3\rt)}{\lz +2 + \sqrt{2}} \rt) \\
&= \frac{\lt(\sqrt{2} + 1\rt) \lt(2+ \sqrt{2} - \sqrt{6}/3\rt)}{\lz +2 + \sqrt{2}}.
\end{align*}
According to Theorem \ref{thmPo}, we have
\begin{align*}
\sup_{\| f \|_{{\rm Lip} (\rz_{1/2})} = 1} M_2 (g_f) &\leq {m_{\rz}} \| \ddz^2 h_{1/2} \|_{\text{\rm Lip}(\rz_{1/2})} + {2} \lt( (2\Xi_2 (\lz)) \wedge \lz^{-1} \rt) \\
&= \frac{\lt(\sqrt{2} + 1\rt) \lt(2+ \sqrt{2} - \sqrt{6}/3\rt)}{\lz +2 + \sqrt{2}} {m_{\rz}} + 2 \lt( (2\Xi_2 (\lz)) \wedge \lz^{-1} \rt). \qed
\end{align*}

\vskip 1cm

\noindent {\bf Proof of Proposition~\ref{example}.}
Let $W_i=W-X_i$, then \cite[(2.27) and (2.29)]{BX06} state that, with $b:=\mu_2$ and $a:=\lambda=\mu-\mu_2$,
\begin{align*}
&\EE\{(f(W-b)-\pz(f))\bone_{W\ge b}\} \notag\\
& =\sum_{i=1}^np_i^2(1-p_i)\EE\{\ddz^2 g_f(W_i-b)\bone_{W_i\ge b}\} \notag\\
 &\ \ \ +g_f(1)\left\{\sum_{i=1}^np_i^2(1-p_i)[\PP(W_i=b-2)-\PP(W_i=b-1)]-a\PP(W=b-1)\right\}\notag\\
& =\sum_{i=1}^np_i^2(1-p_i)\EE\{\ddz^2 g_f(W_i-b)\bone_{W_i\ge b}\} +g_f(1)\EE\{(W-\mu)\bone_{W<b}\},
\end{align*}
which implies
\begin{align}
&\EE\{f(W-b)\bone_{W\ge b} -\pz(f)\}  \notag\\
&=\EE\{(f(W-b)-\pz(f))\bone_{W\ge b}\} -\pz(f)\PP(W< b)\notag\\
& =\sum_{i=1}^np_i^2(1-p_i)\EE\{\ddz^2 g_f(W_i-b)\bone_{W_i\ge b}\} +g_f(1)\EE\{(W-\mu)\bone_{W<b}\}-\pz(f)\PP(W< b).\lb{example4}
\end{align}
Without loss of generality, we assume $f(j)=0$ for all $j\le 0$ so \Ref{stein} ensures $g_f(1)=-\frac1\lambda \pz(f)$ and \Ref{example4} gives
\begin{align}
&\EE\{f((W-b)\bone_{W\ge b})-\pz(f)\} \notag \\
& =\sum_{i=1}^np_i^2(1-p_i)\EE\left\{\ddz^2g_f(W_i-b)\bone_{W_i\ge b}\right\} -\frac{\pz(f)}{\lambda}\EE\{(W-b)\bone_{W<b}\}.\lb{example5}
\end{align}
Using \Ref{rz2}, we have
\begin{align}&\left|\EE\left\{\ddz^2g_f(W_i-b)\bone_{W_i\ge b+1}\right\}\right|\notag\\
&{\leq} \frac{{2}}{1\vee\lambda}\EE\left\{ \ddz\rz(W_i-b)\bone_{W_i\ge b+1}\right\}= \frac{{2}}{1\vee\lambda}\EE\left\{ [2(W_i-b)+1]\bone_{W_i\ge b+1}\right\}\notag\\
&\le {2} +\frac{{4}}{1\vee\lambda}\EE\left\{ (W-b)\bone_{W\ge b+1}\right\}\le {6}+\frac{{4}}{1\vee\lambda}\EE\left\{(b-W)\bone_{W\le b}\right\}.\lb{example6}
\end{align}
On the other hand, \Ref{Mg2} ensures $|\ddz^2g_f(0)|\le2$, which in turn implies
\begin{align}&\left|\sum_{i=1}^np_i^2(1-p_i)\EE\left\{\ddz^2g_f(W_i-b)\bone_{W_i= b}\right\}\right|\notag\\
&\le 2 \sum_{i=1}^np_i^2(1-p_i)\PP(W_i=b)\notag\\
&\le 2 \sum_{i=1}^np_i^2\PP(W=b) \le 2\mu_2 \PP(W\le b).\lb{example7}
\end{align}
Direct verification gives
\begin{align}
\left|\EE\{(W-b)\bone_{W\le b}\}\right|\le \mu_2\PP(W\le b)\lb{example8}
\end{align}
and
\begin{align}
\PP(W\le b)\le e^{-\lambda^2/(2\mu)},\lb{example9}
\end{align}
where the last inequality is due to \cite[Theorem~2.7]{CL06}. The observations that $|f(\cdot)|\le \cdot^2$ implies $|\pz(f)|\le \lambda^2+\lambda$ and $\mu_2-\mu_3 \le \mu- \mz_2 = \lambda$ implies $(\mz_2 - \mz_3) / \lambda \le 1$, and then combining \Ref{example5}, \Ref{example6}, \Ref{example7}, \Ref{example8} and \Ref{example9}, we obtain \Ref{example1}.

For the claim \Ref{example2}, using  \Ref{example9}, we have
$$\mW_2(\law ((W-\mu_2)\bone_{W\ge \mu_2}),\law (W-\mu_2))\le \left\{\EE[(W-\mu_2)^2\bone_{W<\mu_2}\right\}^{1/2}\le \mu_2\PP(W<\mu_2)^{1/2}\le \mu_2e^{-\lambda^2/(4\mu)},$$
hence \Ref{example2} is a direct consequence of the triangle inequality, \Ref{example1} and Proposition~\ref{W2}.
\qed

\section*{Acknowledgements}

Parts of this research were supported by NNSFS of China Nos. 11701588, {11571043, 11431014, 11871008} and ARC Discovery Grant DP150101459.

\bibliographystyle{alpha}

\end{document}